\documentclass[a4paper, 11pt, twoside]{article}
\usepackage{amsmath}
\usepackage{amscd}
\usepackage{amssymb}
\usepackage{amsthm}
\usepackage{graphics}
\usepackage{indentfirst}
\usepackage{latexsym}
\usepackage{amsfonts}
\usepackage{multicol}
\usepackage{stmaryrd}
\usepackage{bbm}
\usepackage{array}
\usepackage{txfonts}

\newtheoremstyle{mattthm}{}{}{\itshape}{}{\bfseries}{.}{ }{}

\theoremstyle{mattthm}
\newtheorem{lemma}{Lemma}[section]
\newtheorem{propn}[lemma]{Proposition}
\newtheorem{thm}[lemma]{Theorem}
\newtheorem{cory}[lemma]{Corollary}

\newtheoremstyle{mattdef}{}{}{}{}{\bfseries}{.}{ }{}

\theoremstyle{mattdef}
\newtheorem*{rmk}{Remark}

\newtheorem{eg}[lemma]{Example}

\newtheorem*{wng}{Warning}

\begin{document}

\newenvironment{pf}{\noindent\textbf{Proof.}}{\hfill \qedsymbol\newline}
\newenvironment{pfof}[1]{\vspace{\topsep}\noindent\textbf{Proof of {#1}.}}{\hfill \qedsymbol\newline}
\newenvironment{pfenum}{\noindent\textbf{Proof.}\indent\begin{enumerate}\vspace{-\topsep}}{\end{enumerate}\vspace{-\topsep}\hfill \qedsymbol\newline}
\newenvironment{pfnb}{\noindent\textbf{Proof.}}{\newline}

\newcommand{\rt}[1]{\rotatebox{90}{$#1$}}

\newcommand{\heit}{\operatorname{ht}}
\newcommand\erange{\{2,3,\dots\}\cup\{\infty\}}
\newcommand{\lo}[1]{\lambda^{(#1)}}
\newcommand\plu{\negthickspace+\negthickspace}
\newcommand\miu{\negthickspace-\negthickspace}
\newcommand\ppmod[1]{\ (\operatorname{mod}\ \ #1)}
\newcommand{\bl}{\boldsymbol{\lambda}}
\newcommand{\bm}{\boldsymbol{\mu}}
\newcommand{\bn}{\boldsymbol{\nu}}
\newcommand{\ol}{\overline}
\newcommand{\ul}{\underline}
\newcommand{\sss}{\mathfrak{S}_}
\newcommand{\mup}{multipartition }
\newcommand{\dom}{\trianglerighteqslant}
\newcommand{\doms}{\vartriangleright}
\newcommand{\ndom}{\ntrianglerighteqslant}
\newcommand{\ndoms}{\not\vartriangleright}
\newcommand{\domby}{\trianglelefteqslant}
\newcommand{\domsby}{\vartriangleleft}
\newcommand{\ndomby}{\ntrianglelefteqslant}
\newcommand{\ndomsby}{\not\vartriangleleft}
\newcommand{\rell}{\sim}
\newcommand{\subs}[1]{\subsection{#1}}
\newcommand{\la}{\lambda}
\newcommand{\fbb}{\mathfrak{B}}
\newcommand{\thmcite}[2]{\textup{\textbf{\cite[#2]{#1}}}\ }
\newcommand{\bbn}{\mathbb{N}}
\newcommand{\bbz}{\mathbb{Z}}
\newcommand{\gs}{\geqslant}
\newcommand{\ls}{\leqslant}
\newcommand{\zez}{\bbz/e\bbz}
\newcommand{\zrz}{\bbz/r\bbz}
\newcommand{\znz}{\bbz/n\bbz}
\newcommand{\zerz}{\bbz/er\bbz}
\newcommand{\hhh}{\mathcal{H}_}
\newcommand{\sect}[1]{\section{#1}}
\newcommand{\ba}{\mathbf{a}}
\newcommand\boa{\mathbf{a}}
\newcommand{\fb}{\mathfrak{b}}
\newcommand{\brls}{\bbz^{r,\ls}}
\newcommand{\bom}{\mathbf{m}}
\newcommand{\bon}{\mathbf{n}}
\newcommand{\bov}{v}
\newcommand\arr\mapsto
\newcommand\rrow\leadsto

\title{Core blocks of Ariki--Koike algebras II: the weight of a core block}
\author{Matthew Fayers\footnote{This research was undertaken with the
support of a Research Fellowship from the Royal Commission for the Exhibition of
1851.  The author is very grateful to the Commission for its generous support.}\\\\Queen Mary, University of London, Mile End Road, London E1 4NS, U.K.\footnote{Correspondence address.}\\\\Massachusetts Institute of Technology, 77 Massachusetts Avenue,\\
Cambridge MA 02139-4307, U.S.A.\footnote{This research was undertaken while the author was visiting M.I.T.\ as a Postdoctoral Fellow.  He is very grateful to Prof.\ Richard Stanley for the invitation, and to M.I.T.\ for its hospitality.}\\\\\texttt{m.fayers@qmul.ac.uk}}
\date{}
\maketitle
\begin{center}
\end{center}
\markboth{Matthew Fayers}{Core blocks of Ariki--Koike algebras II}
\pagestyle{myheadings}

\begin{abstract}
We study combinatorial blocks of multipartitions, exploring further the notions of \emph{weight}, \emph{hub} and \emph{core block} introduced by the author in earlier papers.  We answer the question of which pairs $(w,\theta)$ occur as the weight and hub of a block, and we examine the action of the affine Weyl group on the set of blocks.
\end{abstract}

\sect{Introduction}

The study of the representation theory of the symmetric groups (and, more recently, the Iwahori--Hecke algebras of type $A$) has always been inextricably linked with the combinatorics of partitions.  More recently, the complex reflection group of type $G(r,1,n)$ and its Hecke algebras (the \emph{Ariki--Koike algebras} or \emph{cyclotomic Hecke algebras}) have been studied, and it is clear that there is a similar link to algebraic combinatorics, but with multipartitions playing the r\^ ole of partitions.  This paper is intended as a contribution to the study of the combinatorics of multipartitions, as it relates to the Ariki--Koike algebra.

An important manifestation of multipartition combinatorics is in the block classification for Ariki--Koike algebras.  Given an Ariki--Koike algebra $\hhh n$ of $G(r,1,n)$ and a multipartition $\bl$ of $n$ with $r$ components, there is an important $\hhh n$-module $S^{\bl}$ called a \emph{Specht module}.  Each Specht module lies in one block of $\hhh n$ and each block contains at least one Specht module, so the block classification for $\hhh n$ amounts to deciding when two Specht modules lie in the same block.  Graham and Lehrer \cite{gl} gave a combinatorially-defined equivalence relation $\sim$ on the set of multipartitions, and conjectured that Specht modules $S^{\bl}$ and $S^{\bm}$ lie in the same block of $\hhh n$ if and only if $\bl\sim\bm$.  This was recently proved by Lyle and Mathas \cite{lm}, and their proof uses some of the author's earlier work on the relation $\sim$ and its equivalence classes (which are referred to as \emph{combinatorial blocks} or \emph{blocks of multipartitions}).

Despite this important application, blocks of multipartitions are not well understood, except in the case $r=1$.  In \cite{weight,core}, the author introduced the notions of \emph{weight}, \emph{hub} and \emph{core block} in an attempt to generalise some familiar notions from the case $r=1$.  In this paper, we continue to study these notions, by addressing the question of existence of a block with given weight and hub, and studying the natural action of the affine Weyl group of type $A_{e-1}^{(1)}$ on the set of blocks.

\sect{Background and notation}

\subs{Multipartitions, residues and blocks}\label{mrb}

A \emph{partition} is a non-increasing sequence $\la=(\la_1,\la_2,\dots)$ of non-negative integers such that the sum $|\la|=\la_1+\la_2+\dots$ is finite.  We say that $\la$ is a partition of $|\la|$.  The unique partition of $0$ is referred to as the \emph{empty partition}, and written as $\odot$.

If $\la$ and $\mu$ are partitions, then we say that $\la$ \emph{dominates} $\mu$ and write $\la\dom\mu$ if
\[\la_1+\dots+\la_i\gs\mu_1+\dots+\mu_i\]
for all $i\gs1$.  If $\la$ is a partition, then the \emph{conjugate partition} $\la'$ is given by
\[\la'_j = \left|\left\{i\gs1\ \left|\ \la_i\gs j\right.\right\}\right|.\]
It is easy to show that $\la\dom\mu$ if and only if $\mu'\dom\la'$.

Throughout this paper $r$ is a fixed positive integer.  A \emph{multipartition} is an $r$-tuple $\bl = (\lo 1,\dots,\lo r)$ of partitions, which are called the \emph{components} of $\bl$.  We write $|\bl| = |\lo 1|+\dots+|\lo r|$, and say that $\bl$ is a multipartition of $|\bl|$.  The unique multipartition of $0$ is referred to as the \emph{empty multipartition}, and written as $\varnothing$.

If $\bl$ is a multipartition, the \emph{Young diagram} of $\bl$ is a subset of $\bbn^2\times\{1,\dots,r\}$; we refer to elements of the latter set as \emph{nodes}, and write them in the form $(i,j)_k$, with $i,j\in\bbn$ and $1\ls k\ls r$.  The Young diagram of $\bl$ is the set
\[\left\{(i,j)_k\in\bbn^2\times\{1,\dots,r\}\ \left|\ j\ls \lo k_i\right.\right\},\]
whose elements are called the nodes of $\bl$.  We may abuse notation by not distinguishing between $\bl$ and its Young diagram.  A node $(i,j)_k$ of $\bl$ is \emph{removable} if $\bl\setminus\{(i,j)_k\}$ is the Young diagram of some multipartition, while a node $(i,j)_k$ not in $\bl$ is an \emph{addable node of $\bl$} if $\bl\cup\{(i,j)_k\}$ is the Young diagram of a multipartition.

Now suppose $e\in\{2,3,\dots\}\cup\{\infty\}$ and that $\boa = (a_1,\dots,a_r)\in(\zez)^r$.  If $(i,j)_k$ is a node, then its \emph{residue} is defined to be the element $j-i+a_k$ of $\zez$.  A node of residue $l$ is referred to as an \emph{$l$-node}.  The \emph{content} of a multipartition $\bl$ is the $e$-tuple $(c_l(\bl))_{l\in\zez}$, where $c_l(\bl)$ is the number of nodes of $\bl$ of residue $l$.  If $\bl$ and $\bm$ are multipartitions, then we write $\bl\rell\bm$ if $\bl$ and $\bm$ have the same content.  Clearly $\rell$ is an equivalence relation, and we refer to an equivalence class under this relation as a \emph{block} (of multipartitions).  These definitions depend on the choice of $e$ and $\boa$, and we may write $\bl\rell_{e;\boa}\bm$ or use the terms `$(e;\boa)$-residue' or `$(e;\boa)$-block' if necessary.

\subs{Notational conventions}

It is customary for authors to abuse notation when dealing with the additive group $\znz$, identifying this set with the set of integers $\{0,\dots,n-1\}$, and performing arithmetic modulo $n$.  In this paper we find it more convenient to take a more strict approach: we regard each element of $\znz$ as a set of integers, namely the set $i+n\bbz = \{i+nz\mid z\in\bbz\}$ for some $i\in\bbz$.  The additive group structure on $\znz$ is then derived from addition of sets, defined in the usual way.  We also freely use the additive action of $\bbz$ on $\znz$ without comment.  Note that this interpretation naturally extends to the case $n\ls0$ (though now the cardinality of $\znz=\bbz/(-n)\bbz$ may no longer be equal to $n$).  With this in mind we will allow the multiplication of an element of $\znz$ by any integer $z$, interpreting the result as an element of $\bbz/zn\bbz$; see for example, the definition of the constants $t_i$ in \S\ref{s41}.

In this paper we frequently have cause to consider a variable $e$, which is either an integer greater than or equal to $2$, or is $\infty$.  When $e$ is finite, we may write $\ol{i}$ instead of $i+e\bbz$ for any integer $i$.  When $e=\infty$, the set $\bbz/e\bbz$ should be read simply as $\bbz$.

\subs{`Root space' notation}\label{nota}

In this short section, we outline some notation which will arise later in various guises.  Suppose that $(I,I')$ is one of the following ordered pairs of sets:
\begin{itemize}
\item
$I=I'=\bbz$;
\item
$I=I'=\bbz/n\bbz$, for some positive integer $n$;
\item
$I=\{0,1,\dots,n\}$, $I'=\{1,\dots,n-1\}$, for some positive integer $n$.
\end{itemize}

Write $\bbz^I$ for the free $\bbz$-module with $I$ as a basis; we view an element of $\bbz^I$ as an $I$-tuple $x=(x_i)_{i\in I}$ of integers (of which only finitely many are non-zero, if $I=\bbz$).  For $i\in I'$ we define the element $\alpha_i$ of $\bbz^I$ by
\[(\alpha_i)_j = -\delta_{(i-1)j}+2\delta_{ij}-\delta_{(i+1)j},\]
and set
\[Q^+ = \bigoplus_{i\in I'}\bbz_{\gs0}\alpha_i\]
(where $\bbz_{\gs0}$ denotes the set of non-negative integers).  In the case $I=\bbz$ or $I=\{0,\dots,n\}$, an element $q$ of $Q^+$ can be uniquely written in the form $\sum_{i\in I'}m_i\alpha_i$ with each $m_i\in\bbz_{\gs0}$, and we define $\heit(q) = \sum_{i\in I'}m_i$ in these cases.  In the case where $I=\zez$, we have $\sum_{i\in I}\alpha_i=0$ (and hence $Q^+ = \bigoplus_{i\in I'}\bbz\alpha_i$), so the coefficients $m_i$ in an expression $q=\sum_{i\in I}m_i\alpha_i$ for an element of $Q^+$ are not uniquely determined.  However, there is a unique way to choose the $m_i$ such that all are non-negative and at least one is zero, and we define $\heit(q)=\sum_{i\in I}m_i$ for this particular choice of $m_i$.

Given $x,y\in M$, we write $x\arr y$ if $y-x\in Q^+$.  In order to clarify the setting, we may write
\[x\stackrel{\infty}\arr y,\qquad\phantom{\text{or}}\qquad x\stackrel{n}\arr y,\qquad\text{or}\qquad x\stackrel{0,n}\arr y,\]
as $I$ equals $\bbz$, $\bbz/n\bbz$ or $\{0,\dots,n\}$ respectively.  Notice that $\stackrel{\infty}\arr$ and $\stackrel{0,n}\arr$ are partial orders, while $\stackrel{n}\arr$ is an equivalence relation.

\begin{rmk}
The title of this subsection alludes to the fact that the notation described above resembles notation from the theory of Kac--Moody algebras.  In fact, this phenomenon occurs frequently in the study of blocks of multipartitions.  Indeed, the weight of a multipartition (defined in the next section) can alternatively be written in terms of a standard invariant bilinear form for the Kac--Moody algebra of type $A^{(1)}_{e-1}$.  The implications of this connection to Kac--Moody algebras are still unclear to the author.
\end{rmk}

\subs{Hub and weight}

Suppose $e\in\erange$.  An \emph{$e$-hub} is defined to be an element $\theta$ of the $\bbz$-module $\bbz^{\bbz/e\bbz}$ for which $\sum_{i\in\zez}\theta_i = -r$.  We often just use the term `hub' if the value of $e$ is clear.  A hub $\theta$ is called \emph{negative} if $\theta_i\ls0$ for all $i$.

Now suppose $\boa\in(\zez)^r$, and $\bl$ is a multipartition.  The $(e;\boa)$-hub of $\bl$ is the $e$-tuple $(\theta(\bl))_{i\in\zez}$, where $\theta_i(\bl)$ is defined to be the number of removable $i$-nodes of $\bl$ minus the number of addable $i$-nodes of $\bl$.  We often use the term `hub', if $e$ and $\boa$ are clear.  It is easy to see that the hub of $\bl$ is an $e$-hub, as defined above: $\bl$ has only finitely many nodes and addable nodes, so certainly only finitely many $\theta_i(\bl)$ are non-zero, and it is easy to show that the total number of addable nodes of each component of $\bl$ exceeds the number of removable nodes by $1$, so that $\sum_{i\in\zez}\theta_i(\bl)=-r$.

\begin{wng}
Note that we have made a slight change of notation from \cite{weight,core}; there, the hub of $\bl$ is written as $\delta(\bl)$.  We have made this change to enable us to reserve the symbol $\delta$ for the Kronecker delta.  We hope that no confusion will result.
\end{wng}

We identify a particular hub, which will be very useful later on.  Define $\Theta = \Theta(e;\boa)$ by
\[\Theta_i = -\sum_{j=1}^r\delta_{ia_j}\]
for $i\in\zez$.  That is, $\Theta(e;\boa)$ is the $(e;\boa)$-hub of the empty multipartition.

Now we define the \emph{$(e;\boa)$-weight} (or simply the \emph{weight}) of a multipartition $\bl$.  With the content $(c_i)_{i\in\zez}$ as defined in \S\ref{mrb}, the weight of $\bl$ is defined to be
\[w(\bl) = \sum_{j=1}^rc_{a_j}(\bl) - \frac12\sum_{i\in\zez}(c_i(\bl)-c_{i-1}(\bl))^2.\]
Though it is not obvious from this formula, $w(\bl)$ is always a non-negative integer.  What is clear from the formula is that any two partitions in the same block have the same weight.  In fact, it is proved in \cite[Proposition 3.2]{weight} that two multipartitions lie in the same block if and only if they have the same hub and the same weight.  So we may speak of the weight and hub of a block $B$ (meaning the weight and hub of any multipartition lying in $B$), which we write as $w(B)$ and $\theta(B)$.

The main question that we address in this paper is: for which hubs $\theta$ and weights $w$ does there exist a block with weight $w$ and hub $\theta$?

\subs{Beta-numbers}\label{betacore}

In this section we fix $e\in\erange$ and $\boa\in(\zez)^r$.  A \emph{multi-charge} (for $\boa$) is defined to be an $r$-tuple $\hat{\boa} = (\hat{a}_1,\dots,\hat{a}_r)\in\bbz^r$ such that for each $i\in\zez$ we have $a_i=\ol{\hat{a}_i}$.  (If $e=\infty$, then there is only one possible choice of multi-charge, namely $\hat{\boa}=\boa$.)  Given a multi-charge $\hat{\boa}$ and a multipartition $\bl$, we define the \emph{beta-numbers}
\[\beta^j_i =\lo j_i+\hat{a}_j-i\]
for all $i\gs 1$ and all $j\in\{1,\dots,r\}$.  For any $j$, the beta-numbers $\beta_1^j,\beta_2^j,\dots$ are distinct, and we refer to the sets
\[B^j = \left\{\left.\beta^j_i\ \right|\ i\gs1\right\}\]
as the \emph{beta-sets} for $\bl$ corresponding to $\hat\boa$.  Note that for each $j$, $B^j$ has the following property.
\[\text{For $i$ sufficiently large, the number of elements of $B^j$ greater than or equal to $-i$ is $\hat a_j+i$.}\tag*{($\ast$)}\]
Conversely, it is easy to see that given a multicharge $\hat\boa$, any $r$ subsets $B^1,\dots,B^r$ of $\bbz$ satisfying ($\ast$) for each $j$ are the beta-sets of some multipartition.

There is an important link between beta-numbers and addable and removable nodes.  Given beta-sets $B^1,\dots,B^r$ for a multipartition $\bl$, it is easy to see that there is a removable $l$-node in row $i$ of the $j$th component of $\bl$ if and only if $\ol{\beta^j_i} = l$ and $\beta^j_i-1\notin B^j$.  Furthermore, if $\bm$ is the multipartition obtained by removing this removable node, then the beta-sets for $\bm$ are
\[B^1,\dots,B^{j-1},B^j\setminus\{\beta^j_i\}\cup\{\beta^j_i-1\},B^{j+1},\dots,B^r.\]
This enables us to express the hub of a multipartition in terms of its beta-numbers.

\begin{lemma}
Suppose $\bl$ is a multipartition, with beta-sets $B^1,\dots,B^r$.  Then for $l\in\zez$ we have
\[\theta_l(\bl) = \sum_{j=1}^r\left|\left\{\beta\in B^j\ \left|\ \ol{\beta}= l,\ \beta-1\notin B^j\right.\right\}\right|-\left|\left\{\beta\in B^j\ \left|\ \ol{\beta}=l-1,\ \beta+1\notin B^j\right.\right\}\right|.\]
\end{lemma}

\begin{cory}\label{remalpha}
Suppose $\bl$ and $\bm$ are multipartitions, and $\bm$ is obtained by removing a removable $l$-node from $\bl$.  Then for $i\in\zez$ we have
\[\theta_i(\bm) = \theta_i(\bl)+\delta_{(i-1)l}-2\delta_{il}+\delta_{(i+1)l}.\]
\end{cory}

\subs{Multicores and core blocks}

Suppose $e$ and $\boa$ are as above, and $\bl$ is a multipartition, with beta-sets $B^1,\dots,B^r$.  If $e$ is finite, we say that $\bl$ is a \emph{multicore} if there do not exist $i\in\bbz$ and $j\in\{1,\dots,r\}$ such that $i\in B^j$ and $i-e\notin B^j$.  It is easy to show that this condition is independent of the choice of multi-charge $\hat{\boa}$.  In fact, $\bl$ is a multicore if and only if each component $\la^{(i)}$ of $\bl$ is an $e$-core, i.e.\ none of the hook lengths of $\la^{(i)}$ is divisible by $e$.  If $e=\infty$, then we deem every multipartition to be a multicore.

We define some more notation, which arises from examining beta-numbers modulo $e$: if $e$ is finite, then for each $i\in\zez$ and each $j\in\{1,\dots,r\}$ we write
\begin{align*}
\fb_{ij}^{\hat{\boa}}(\bl) &= \max\left\{\left.\beta\in B^j\ \right|\ \ol{\beta}=i\right\}.\\
\intertext{If $e=\infty$, then we define}
\fbb_{ij}(\bl) &= \begin{cases}
1 & (i\in B^j)\\
0 & (i\notin B^j)
\end{cases}\end{align*}
for $i\in\bbz$, $j\in\{1,\dots,r\}$.  We can now recall one of the main results from \cite{core}, which motivates the definition of a core block.

\begin{thm}\thmcite{core}{Theorem 3.1}\label{coremain}
Suppose $e$ is finite, and that $\bl$ is a multipartition lying in a block $B$.  Then the following are equivalent.
\begin{enumerate}
\item
$\bl$ is a multicore, and there exist a multi-charge $\hat{\boa}$ and integers $(b_i)_{i\in\zez}$ such that for each $i\in\zez$ and $j\in\{1,\dots,r\}$ we have
\[\fb^{\hat{\boa}}_{ij}(\bl) \in \{b_i,b_i+e\}.\]
\item\label{nolower}
There is no block with the same hub as $B$ and smaller weight.
\item\label{everycore}
Every multipartition in $B$ is a multicore.
\end{enumerate}
\end{thm}

If $e<\infty$, we say that a block is a \emph{core block} if the conditions of Theorem \ref{coremain} are satisfied for any (and hence every) $\bl\in B$.  If $e=\infty$, then every block is deemed to be a core block; it follows from \cite[Proposition 1.3]{core} that when $e=\infty$ a block is uniquely determined by its hub, so that property (\ref{nolower}) of Theorem \ref{coremain} holds for every block when $e=\infty$; property (\ref{everycore}) holds too, since every multipartition is defined to be a multicore.

The idea of the definition is that a core block should behave like `a block at $e=\infty$'.  In fact, this can be made precise, via the notion of a `lift' of a core block (see \cite[\S3.4]{core}).  When $r=1$, a block $B$ is a core block if and only if it has weight $0$, which happens if and only if $|B|=1$; furthermore, every multicore lies in a core block.

\begin{eg}
Suppose $r=3$ and $e=4$, and $\boa=(\ol0,\ol0,\ol2)$.  Let $\bl=((2),\odot,(1))$.  Given the multi-charge $\hat\boa = (0,0,2)$, we get the beta-sets
\begin{align*}
B^1 &= \{\dots,-3,-2,1\},\\
B^2 &= \{\dots,-3,-2,-1\},\\
B^3 &= \{\dots,-3,-2,-1,0,2\}.
\end{align*}
Hence
\begin{alignat*}3
\fb_{\ol01}^{\hat{\boa}}(\bl) &= -4,&\qquad \fb_{\ol02}^{\hat{\boa}}(\bl) &= -4,&\qquad \fb_{\ol03}^{\hat{\boa}}(\bl) &= 0,\\
\fb_{\ol11}^{\hat{\boa}}(\bl) &= 1,&\qquad \fb_{\ol12}^{\hat{\boa}}(\bl) &= -3,&\qquad \fb_{\ol13}^{\hat{\boa}}(\bl) &= -3,\\
\fb_{\ol21}^{\hat{\boa}}(\bl) &= -2,&\qquad \fb_{\ol22}^{\hat{\boa}}(\bl) &= -2,&\qquad \fb_{\ol23}^{\hat{\boa}}(\bl) &= 2,\\
\fb_{\ol31}^{\hat{\boa}}(\bl) &= -5,&\qquad \fb_{\ol32}^{\hat{\boa}}(\bl) &= -1,&\qquad \fb_{\ol33}^{\hat{\boa}}(\bl) &= -1,\\
\end{alignat*}
and we see that condition (1) in Theorem \ref{coremain} holds, so $\bl$ lies in a core block.  It is easy to see that (3) holds too: any multipartition $\bm$ in the same block as $\bl$ has $|\bm|=|\bl|=3$, and it is easy to show in general that if $|\bm|<e$ then $\bm$ is a multicore.
\end{eg}

The use of core blocks will help us to answer our main question.  The fact that the core blocks are the blocks satisfying property (\ref{nolower}) of Theorem \ref{coremain} implies that if there exists a block $B$ with hub $\theta$ then there exists a core block with hub $\theta$; since a block is determined by its hub and weight, this core block is unique.  So we may speak of \emph{the core block with hub $\theta$}, when we know that there is some block with hub $\theta$.  If $B$ is a block with hub $\theta$ and $C$ is the core block with this hub, then it follows from \cite[Proposition 1.3]{core} that $w(B)-w(C)$ is an integer multiple of $r$.  Conversely, given a core block $C$ and a non-negative integer $n$, it is easy to construct a block with the same hub and weight $w(C)+nr$ (by adding $n$ rim $e$-hooks to a multipartition in $C$), so we have the following.

\begin{lemma}\label{listwt}
Suppose $e<\infty$.  Then there exists a block with hub $\theta$ and weight $w$ if and only if
\begin{itemize}
\item
there exists a core block $C$ with hub $\theta$, and
\item
$w-w(C)$ is a non-negative integer multiple of $r$.
\end{itemize}
\end{lemma}

In view of this, we can re-state our main question as follows.  Given a hub $\theta$, does there exist a block with hub $\theta$, and if so, what is the weight of the core block with hub $\theta$?

\begin{rmk}
Blocks are actually slightly peripheral to this paper, in the sense that it would be possible to state all the results in terms of multipartitions; in \cite{lm}, Lyle and Mathas introduced the term `reduced multicore' for a multipartition lying in a core block, and so we could ask: given a hub $\theta$, does there exist a multipartition with hub $\theta$, and if so, what is the weight of a reduced multicore with hub $\theta$?  However, since the notions of hub and weight were introduced precisely in order to study blocks, the results here are naturally stated in terms of blocks.
\end{rmk}

We conclude this section by observing that the hub of a multipartition $\bl$ can be recovered from the integers $\fb^{\hat{\boa}}_{ij}(\bl)$ or $\fbb_{ij}(\bl)$.

\begin{lemma}\thmcite{core}{Lemma 3.2 \& proof of Proposition 3.6}\label{recoverhub}
Suppose $e$ and $\boa$ are as above, and $\hat{\boa}$ is a multi-charge for $\boa$.
\begin{enumerate}
\item
If $e$ is finite, then
\[\theta_i(\bl) = \frac{\sum_{j=1}^r\left(\fb^\ba_{ij}(\bl)-\fb^\ba_{(i-1)j}(\bl)\right)-r}e\]
for $i\in\zez$.
\item
If $e=\infty$, then
\[\theta_i(\bl) = \sum_{j=1}^r\left(\fbb_{ij}(\bl)-\fbb_{(i-1)j}(\bl)\right)\]
for $i\in\bbz$.
\end{enumerate}
\end{lemma}

%
%
%
%

\subs{Actions of the affine Weyl group}

Given $e$ as above, let $W_e$ denote the Weyl group of type $A^{(1)}_{e-1}$ (or type $A_\infty$, if $e=\infty$).  This is the Coxeter group with generators $s_i$ for $i\in\zez$, and relations
\begin{alignat*}2
s_i^2 &=1&\qquad&\text{for all $i$},\\
s_is_j &= s_js_i&\qquad &\text{whenever $i\neq j\pm1$},\\
s_is_js_i &= s_js_is_j&\qquad&\text{whenever $i=j+1\neq j-1$}.
\end{alignat*}
There are several actions of $W_e$ which are of interest to us.  First we recall the well-known permutation action on $\bbz$, which gives $W_e$ the name `generalised symmetric group' (if $e<\infty$) or `finitary symmetric group' (if $e=\infty$).

\begin{lemma}\label{zaction}
There is an action of $W_e$ on $\bbz$ given by
\[s_i(j) = \begin{cases}
j+1 & (\ol j=i-1)\\
j-1 & (\ol j=i)\\
j & (\text{otherwise})\text,
\end{cases}\]
for all $i\in\zez$ and $j\in\bbz$.
\end{lemma}

Next we give an action on hubs.  The following lemma is straightforward; the slightly nebulous expression using the Kronecker delta obviates the need for a separate definition for the case $e=2$.

\begin{lemma}\label{hubaction}
There is an action of $W_e$ on the set of $e$-hubs, given by
\[\left(s_i(\theta)\right)_j = \theta_j + \theta_i\left(\delta_{(i-1)j}-2\delta_{ij}+\delta_{(i+1)j}\right)\text,\]
for a hub $\theta$ and $i,j\in\zez$.
\end{lemma}

Now we describe an action on multipartitions, which yields an action on blocks.  Fix $e\in\erange$ and $\boa\in(\zez)^r$.

\begin{propn}\label{mpaction}
Given a multipartition $\bl$ and $i\in\zez$, define $s_i(\bl)$ to be the multipartition obtained by simultaneously removing all removable $i$-nodes and adding all addable $i$-nodes.  This defines an action of $W_e$ on the set of multipartitions.  Moreover, we have
\[w(s_i(\bl)) = w(\bl),\qquad\theta(s_i(\bl)) = s_i(\theta(\bl))\]
for any multipartition $\bl$, where $s_i(\theta(\bl))$ is as defined in Lemma \ref{hubaction}.
\end{propn}

\begin{pf}
Choose a multicharge $\hat{\boa}$, and for $j\in\{1,\dots,r\}$ define the beta-sets $B^j(\bl)$ and $B^j(s_i(\bl))$ for $\bl$ and $s_i(\bl)$ as in \S\ref{betacore}.  As noted in \cite[\S4.2.1]{weight}, we have
\[B^j(s_i(\bl)) = s_i\left(B^j(\bl)\right),\]
where the action on the right-hand side is the permutation action of $W_e$ on $\bbz$ described in Lemma \ref{zaction}.  This implies that the definition of $s_i(\bl)$ gives an action.  The statements concerning the weight and hub of $s_i(\bl)$ are proved in \cite[Proposition 4.6]{weight}.
\end{pf}

Since two partitions lie in the same block if and only if they have the same hub and weight, we see that the action of $W_e$ on multipartitions given in Proposition \ref{mpaction} reduces to an action on the set of blocks, with the property that
\[w(s_i(B))=w(B)\]
and
\[\theta(s_i(B)) = s_i(\theta(B))\]
for any block $B$ and any $i\in\zez$.

We give a corollary which is of great importance to the main question in this paper.

\begin{cory}\label{existorb}
Suppose $e,\boa$ are chosen as above, and that $\theta$ and $\kappa$ are hubs lying in the same $W_e$-orbit.  Then: \begin{enumerate}
\item\label{eo1}
there exists a block with hub $\theta$ if and only if there exists a block with hub $\kappa$;
\item\label{eo2}
if there exist such blocks, then the core blocks with hubs $\theta$ and $\kappa$ lie in the same $W_e$-orbit and have the same weight.
\end{enumerate}
\end{cory}

\begin{pf}
It suffices to assume that $\kappa=s_i(\theta)$ for some $i\in\zez$.
\begin{enumerate}
\item
If $B$ is a block with hub $\theta$, then (from above) the block $s_i(B)$ has hub $\kappa$.  The `if' part is similar.
\item
Let $B$ and $C$ be the core blocks with hubs $\theta$, $\kappa$ respectively.  The block $s_i(B)$ has hub $\kappa$, and since $C$ has the smallest weight of any block with hub $\kappa$, we have
\[w(C)\ls w(s_i(B)) = w(B).\]
Similarly, we have $w(B)\ls w(C)$, so that $w(B)=w(C)$.  This gives $w(C)=w(s_i(B))$, and since a block is determined by its hub and weight, we deduce that $C=s_i(B)$.
\end{enumerate}
\end{pf}

\sect{The case $e=\infty$}\label{sec3}

Throughout this section we fix $e=\infty$ and $\boa\in\bbz^r$, and answer our main question in this case.

\subs{Negative hubs}

We address first the existence of blocks with negative hubs; it will then be a simple matter to extend to the general case, using the action of $W_\infty$.

Using the integers $\fbb_{ij}(\bl)$, we can interpret the question of the existence of a multipartition (and hence a block) with hub $\theta$ by looking at the existence of zero--one matrices.  Suppose that $\theta$ is a negative hub, and let $d_i = d_i(\theta)$ be defined for $i\in\bbz$ by
\[d_i = r+\sum_{l=-\infty}^i\theta_l.\]
Then the sequence $(d_i)_{i\in\bbz}$ is weakly decreasing, and for sufficiently large $i$ we have $d_{-i}=r$ and $d_i=0$.  Choose and fix $k$ sufficiently small that $d_k=r$ and $a_j\gs k$ for each $j$, and $l$ sufficiently large that $d_l=0$ and $a_j\ls l$ for each $j$.

If there is a multipartition $\bl$ with hub $\theta$, then by Lemma \ref{recoverhub} we have
\[d_i = \sum_{j=1}^r\fbb_{ij}(\bl)\]
for all $i$.    In particular, $\fbb_{ij}(\bl)=1$ for all $j$ when $i<k$, and $\fbb_{ij}(\bl)=0$ for all $j$ when $i>l$.  So by property ($\ast$) from \S\ref{betacore} the sum
\[\fbb_{kj}(\bl)+\fbb_{(k+1)j}(\bl)+\dots+\fbb_{lj}(\bl)\]
equals $a_j-k$, for each $j\in\{1,\dots,r\}$.  Hence the $(l-k+1)\times r$ matrix with $(i,j)$-entry $\fbb_{(i+k-1)j}(\bl)$ is a zero--one matrix with row sums $d_k\gs \dots\gs d_l$ from top to bottom, and columns sums $a_1-k,\dots,a_r-k$ from left to right.

Conversely, suppose we have a zero--one matrix $A$ with these row and column sums.  Then it easy to construct a multipartition $\bl$ with hub $\theta$, by setting
\[\fbb_{ij}(\bl) = \begin{cases}
1 & (i<k)\\
A_{(i+1-k)j} & (k\ls i\ls l)\\
0 & (i>l)
\end{cases}\]
for each $i,j$.  So we have shown the following.

\begin{lemma}
\label{mpmx}
Suppose $\theta$ is a negative hub, and choose $k,l$ as above.  There exists a block with hub $\theta$ if and only if there exists an $(l-k+1)\times r$ zero--one matrix with columns sums $a_1-k,\dots,a_r-k$ from left to right, and row sums $d_k,\dots,d_l$ from top to bottom.
\end{lemma}

Necessary and sufficient conditions for the existence of zero--one matrices with prescribed row and column sums are given by the Gale--Ryser Theorem \cite{gale,ryser}.  For each integer $t\in\{k,\dots,l\}$, let $c_t$ be the number of values $j$ for which $a_j>t$.  The Gale--Ryser Theorem can be stated in our notation as follows.

\begin{thm}\label{gr}
A zero--one matrix as in Lemma \ref{mpmx} exists if and only if
\[c_k+\dots+c_t\gs d_k+\dots+d_t\]
for each $t\in\{k,\dots,l\}$, with equality when $t=l$.
\end{thm}

In other words, if we let $\sigma$ be the partition $(c_k,c_{k+1},\dots,c_l,0,0,\dots)$, and $\tau$ the partition $(d_k,d_{k+1},\dots)$, then there exists a block with hub $\theta$ if and only if $|\sigma|=|\tau|$ and $\sigma\dom\tau$.  Since the dominance order is reversed under conjugation of partitions, this is equivalent to saying that $|\sigma'|=|\tau'|$ and $\sigma'\domby\tau'$.  The latter interpretation will be more convenient for our purposes.

We wish to re-state the conditions of Theorem \ref{gr} in a way which is more convenient for generalising to arbitrary hubs.  To begin with, re-arrange $a_1,\dots,a_r$ in ascending order as $b_1\ls\dots\ls b_r$.  Since $\theta$ is a negative hub, there exist integers $f_1\ls\dots\ls f_r$ such that
\[\theta_i = -\left|\left\{j\in\{1,\dots,r\}\ \left|\ f_j=i\right.\right\}\right|\]
for each $i$.  (If $\bl$ is a multipartition with hub $\theta$, then $f_1,\dots,f_r$ may be interpreted as the residues of the addable nodes of $\bl$ with multiplicity, where `with multiplicity' entails cancelling addable nodes with removable nodes.)  Note that we have
\[d_i = \left|\left\{j\in\{1,\dots,r\}\ \left|\ f_j>i\right.\right\}\right|\]
for each $i$, which gives the following.

\begin{propn}
\label{mainneg}
Suppose $\theta$ is a negative hub, and let $b_1,\dots,b_r,f_1,\dots,f_r$ be as above.  Then there exists a block with hub $\theta$ if and only if
\[b_1+\dots+b_j\gs f_1+\dots+f_j\]
for each $j\in\{1,\dots,r\}$, with equality when $j=r$.
\end{propn}

\begin{pf}
We have
\begin{align*}
\sigma' &= (b_r-k,b_{r-1}-k,\dots,b_1-k,0,0,\dots)\\
\intertext{and}
\tau' &= (f_r-k,f_{r-1}-k,\dots,f_1-k,0,0,\dots).
\end{align*}
The condition $|\sigma'|=|\tau'|$ is therefore equivalent to
\[b_1+\dots+b_r = f_1+\dots+f_r,\]
and  the condition $\sigma'\domby\tau'$ is equivalent to
\[b_r+b_{r-1}+\dots+b_t\ls f_r+f_{r-1}+\dots+f_t\]
for all $t$, which (in the presence of the condition $|\sigma'|=|\tau'|$) is equivalent to
\[b_1+\dots+b_j\gs f_1+\dots+f_j\]
for each $j\in\{1,\dots,r\}$.
\end{pf}

We wish to generalise this to arbitrary hubs, and give an expression for the weight of a block in terms of its hub.  
First we give a different statement of the conditions in Proposition \ref{mainneg}, which will be useful in the next section.  This requires some more notation.

Write $\brls$ for the set of $r$-tuples $\bom = (m_1,\dots,m_r)$ of integers such that $m_1\ls\dots\ls m_r$.  Given $\bom,\bon\in\brls$, write $\bom\rrow\bon$ if there exist $1\ls p<q\ls r$ such that
\[n_j = \begin{cases}
m_p-1 & (j=p)\\
m_q+1 & (j=q)\\
m_j & (\text{otherwise})
\end{cases}\]
for $j=1,\dots,r$.  Note that this condition automatically implies that $m_{p-1}<m_p$ if $p>1$, and $m_q<m_{q+1}$ if $q<r$.  Informally, $\bom\rrow\bon$ if $\bon$ is obtained from $\bom$ by moving two values $m_p,m_q$ apart without changing the order of $m_1,\dots,m_r$.

Now recall the definition of $\Theta=\Theta(e;\boa)$, and the notation $x\stackrel{\infty}\arr y$ from \S\ref{nota}.

\begin{propn}\label{tfae}
Suppose $\theta$ is a negative hub, and let $b_1\ls\dots\ls b_r$ and $f_1\ls\dots\ls f_r$ be as above.  Then the following are equivalent.
\begin{enumerate}
\item
There exists a block with hub $\theta$.
\item
We have $b_1+\dots+b_j\gs f_1+\dots+f_j$ for any $j\in\{1,\dots,r\}$, with equality when $j=r$.
\item
There exist $\bom_0,\dots,\bom_s\in\bbz^{r,\ls}$ such that
\[(b_1,\dots,b_r)=\bom_0\rrow\dots\rrow\bom_s=(f_1,\dots,f_r).\]
\item
$\Theta\stackrel{\infty}\arr\theta$.
\end{enumerate}
\end{propn}

\begin{pf}
\begin{description}
\item[\textnormal{(2)$\Rightarrow$(1)}]
This is part of Proposition \ref{mainneg}.
\item[\textnormal{(3)$\Rightarrow$(2)}]
Since the relation in (2) is transitive, we may assume that $(b_1,\dots,b_r)\rrow(f_1,\dots,f_r)$.  So suppose $1\ls p<q\ls r$ and
\[f_j = \begin{cases}
b_p-1 & (j=p)\\
b_q+1 & (j=q)\\
b_j & (\text{otherwise}).
\end{cases}\]
Then for $1\ls j\ls r$, we have
\[f_1+\dots+f_j = \begin{cases}
b_1+\dots+b_j & (j<p)\\
b_1+\dots+b_j-1 & (p\ls j<q)\\
b_1+\dots+b_j & (q\ls j),
\end{cases}\]
and (2) holds.
\item[\textnormal{(4)$\Rightarrow$(3)}]
Suppose $\theta = \Theta+\sum_{i\in\bbz}m_i\alpha_i$, with each $m_i$ non-negative.  If every $m_i$ equals $0$ then the result is trivial, so we assume otherwise.  Let $s$ be minimal such that $m_s>0$, and let $t\gs s$ be minimal such that $m_{t+1}=0$.  Set $\kappa = \theta-(\alpha_s+\dots+\alpha_t)$.  Then we claim that $\kappa$ is a negative hub.  Since $\theta$ is a negative hub, we just need to show that $\kappa_{s-1},\kappa_{t+1}\ls0$, i.e.\ that $\theta_{s-1}$ and $\theta_{t+1}$ are strictly negative.  But since $m_{s-1} = m_{t+1}=0$ while $m_s,m_t>0$, the $(s-1)$th and $(t+1)$th entries of $\sum_{i\in\bbz}m_i\alpha_i$ are strictly negative, and the fact that $\Theta_{s-1},\Theta_{t+1}\ls0$ completes the proof of the claim.

So $\kappa$ is a negative hub, and $\Theta\arr\kappa\arr\theta$.  If we let $g_1\ls \dots\ls g_r$ be the integers such that $\kappa_i = -\left|\{j\in\{1,\dots,r\}\mid g_j=i\}\right|$, then by induction on $\sum_i m_i$, condition (3) holds with $g_1,\dots,g_r$ in place of $f_1,\dots,f_r$.  So, by the transitivity of the relation in (3), it suffices to show that $(g_1,\dots,g_r)\rrow(f_1,\dots,f_r)$.

Note that we have
\[\theta_i = \kappa_i+\delta_{is}-\delta_{i(s-1)}+\delta_{it}-\delta_{i(t+1)}\]
for each $i$; in particular, $\kappa_s$ and $\kappa_t$ are strictly negative, i.e.\ there exist $p$ and $q$ such that $g_p=s$ and $g_q=t$.  Letting $p$ be minimal such that $g_p=s$, and $q$ maximal such that $g_q=t$, we set
\[f_j' = \begin{cases}
g_p-1 & (j=p)\\
g_q+1 & (j=q)\\
g_j & (\text{otherwise}).
\end{cases}\]
Then we have $f_1'\ls\dots\ls f_r'$ and $\theta_l = -\left|\left\{j\in\{1,\dots,r\}\ \left|\ f_j'=l\right.\right\}\right|$, which means that $f_j'=f_j$ for all $j$ and we are done.
\item[\textnormal{(1)$\Rightarrow$(4)}]
We show that (1)$\Rightarrow$(4) for any hub $\theta$, not just a negative one.  If there exists a multipartition $\bl$ with hub $\theta$, then we can reach the empty multipartition from $\bl$ by successively removing removable nodes.  Each time we remove an $i$-node, we subtract $\alpha_i$ from the hub, by Corollary \ref{remalpha}.  Since the empty multipartition has hub $\Theta$, the result follows.
\end{description}
\end{pf}

\begin{eg}\label{ex1}
Suppose $r=4$ and $\boa=(1,2,1,0)$.  Define $\theta$ by
\[\theta_i = \begin{cases}
-2 & (i=0)\\
-1 & (i=1)\\
-1 & (i=3)\\
0 & (\text{otherwise}).
\end{cases}\]
Then $\theta$ is a negative hub, and we may verify the conditions of Proposition \ref{tfae} for $\theta$.  First we note that the multipartition $(\odot,(1),(1),\odot)$ has hub $\theta$: its Young diagram, with the residue of nodes and addable nodes marked, may be drawn as follows.
\[
\begin{picture}(30,30)
\put(0,20){\line(0,1){10}}
\put(0,30){\line(1,0){10}}
\put(5,19){$1$}
\end{picture}\qquad
\begin{picture}(30,30)
\put(0,5){\line(0,1){25}}
\put(15,15){\line(0,1){15}}
\put(0,15){\line(1,0){15}}
\put(0,30){\line(1,0){25}}
\put(5,4){$1$}
\put(5,19){$2$}
\put(20,19){$3$}
\end{picture}\qquad
\begin{picture}(30,30)
\put(0,5){\line(0,1){25}}
\put(15,15){\line(0,1){15}}
\put(0,15){\line(1,0){15}}
\put(0,30){\line(1,0){25}}
\put(5,4){$0$}
\put(5,19){$1$}
\put(20,19){$2$}
\end{picture}\qquad
\begin{picture}(30,30)
\put(0,20){\line(0,1){10}}
\put(0,30){\line(1,0){10}}
\put(5,19){$0$}
\end{picture}
\]
So (1) holds.  For (2) and (3), we note that $(b_1,b_2,b_3,b_4)=(0,1,1,2)$ and $(f_1,f_2,f_3,f_4)=(0,0,1,3)$, so (2) and (3) are readily verified.  For (4), we have
\[(\theta-\Theta)_i = \begin{cases}
-1 & (i=0)\\
1 & (i=1)\\
1 & (i=2)\\
-1 & (i=3)\\
0 & (\text{otherwise}),
\end{cases}\]
so that $\theta = \Theta+\alpha_1+\alpha_2$.
\end{eg}

\subs{Arbitrary hubs}

Now we give necessary and sufficient conditions on an arbitrary hub $\theta$ for the existence of a block with hub $\theta$, and give an expression for the weight of such a block.

As before, we define
\[d_i=d_i(\theta) = r+\sum_{l=-\infty}^i\theta_l\]
for $i\in\bbz$.  And as before, if $\bl$ is a multipartition with hub $\theta$, then $d_i$ will be the sum $\fbb_{i1}(\bl)+\dots+\fbb_{ir}(\bl)$.  In particular, for such a $\bl$ to exist we must have $0\ls d_i\ls r$ for each $i$; we say that $\theta$ is \emph{plausible} if this condition holds.  
Now for $j=1,\dots,r-1$ define
\begin{alignat*}3
v_j(\theta) &= |\{i\in\bbz\mid\ &d_i=j\}&|.&&\\
\intertext{Also define}
v_0(\theta) &= |\{i<0\mid &d_i=0\}&|-|\{i\gs0\mid d_i\neq0&&\}|,\\
v_r(\theta) &= |\{l\gs 0\mid &d_i=r\}&|-|\{i<0\mid d_i\neq r&&\}|.
\end{alignat*}
Note that since we have $d_i=0$ and $d_{-i}=r$ for $i$ sufficiently large, $v_0(\theta),v_1(\theta),\dots,v_r(\theta)$ are well-defined integers.  Moreover, we have $v_0(\theta)+\dots+v_r(\theta)=0$ if and only if $\theta$ is plausible.  We let $\bov(\theta)=(v_0(\theta),\dots,v_r(\theta))$, regarded as an element of $\bbz^{\{0,\dots,r\}}$.

Recalling from Lemma \ref{hubaction} the $W_\infty$-action on the set of hubs, we have the following.

\begin{lemma}\label{actv}
The action of $W_\infty$ on the set of hubs preserves the set of plausible hubs.  If $\theta,\kappa$ are plausible hubs, then $\theta$ and $\kappa$ lie in the same $W_\infty$-orbit if and only if $\bov(\theta)=\bov(\kappa)$.  Each $W_{\infty}$-orbit of plausible hubs contains exactly one negative hub.
\end{lemma}

\begin{pf}
If $i,l\in\bbz$, then it is easy to check that
\[d_l(s_i\theta) = \begin{cases}
d_{i-1}(\theta) & (l=i)\\
d_i(\theta) & (l=i-1)\\
d_l(\theta) & (\text{otherwise}).
\end{cases}\]
Hence two hubs $\theta$ and $\kappa$ lie in the same $W_\infty$-orbit if and only if there is some finitary permutation $\pi$ of $\bbz$ such that $d_l(\kappa) = d_{\pi(l)}(\theta)$ for all $l$; this proves the first two statements.  For the third statement, observe that $\theta$ is negative if and only if the sequence $(d_i\theta)$ is weakly decreasing.
\end{pf}

Now we can give our main result for $e=\infty$; recall the notation $x\stackrel{0,r}\arr y$ and $\heit(y-x)$ from \S\ref{nota}.

\begin{thm}\label{maininf}
Suppose $\theta$ is a hub.  Then there exists a block with hub $\theta$ if and only if $\bov(\Theta)\stackrel{0,r}\arr\bov(\theta)$, and if this is the case then the unique such block has weight $\heit(\bov(\theta)-\bov(\Theta))$.
\end{thm}

\begin{pf}
Let us suppose first that $\bov(\Theta)\arr\bov(\theta)$.  This implies that 
\[v_0(\theta)+\dots+v_r(\theta) = v_0(\Theta)+\dots+v_r(\Theta);\]
since $\Theta$ is obviously plausible, we see that $\theta$ is plausible.  In proving the existence of a block with hub $\theta$, we may (in view of Lemma \ref{actv}) replace $\theta$ with any other hub in the same $W_\infty$-orbit; in particular, since any orbit of plausible hubs contains a negative hub, we may assume that $\theta$ is negative.  Given this assumption, we will show that there exists a block with hub $\theta$ by verifying condition (2) of Proposition \ref{tfae}.

Letting $b_1\ls\dots\ls b_r$ and $f_1\ls\dots\ls f_r$ be as in the previous section, we have
\[v_j(\Theta) = b_{r+1-j}-b_{r-j},\qquad v_j(\theta) = f_{r+1-j}-f_{r-j}\]
for $j=0,\dots,r$, where we interpret $b_0,b_{r+1},f_0,f_{r+1}$ as zero.  Hence we have
\begin{alignat*}4
v_{r+1-j}(\Theta)&+2v_{r+2-j}(\Theta)&&+\dots+jv_r(\Theta) &&= & b_1&+\dots+b_j,\\
v_{r+1-j}(\theta)&+2v_{r+2-j}(\theta)&&+\dots+jv_r(\theta) &&= & f_1&+\dots+f_j
\end{alignat*}
for each $j$.  The condition $\bov(\Theta)\arr\bov(\theta)$ says that we have $\bov(\theta) = \bov(\Theta)+\sum_{j=1}^{r-1}m_j\alpha_j$ with each $m_j$ non-negative, and we find
\begin{align*}
(b_1+\dots+b_j)-(f_1+\dots+f_j) &= \left(v_{r+1-j}(\Theta)-v_{r+1-j}(\theta)\right)+2\left(v_{r+2-j}(\Theta)-v_{r+2-j}(\theta)\right)+\dots+j\Big(v_r(\Theta)-v_r(\theta)\Big)\\
&= \begin{cases}
m_{r-j} & (1\ls j\ls r-1)\\
0 & (j=r).
\end{cases}
\end{align*}
So we have
\[b_1+\dots+b_j\gs f_1+\dots+f_j\]
for all $j$, with equality when $j=r$, and so by Proposition \ref{mainneg} there exists a block with hub $\theta$.

Conversely, suppose $\bl$ is a multipartition with hub $\theta$.  We will prove by induction on $|\bl|$ that $\bov(\Theta)\arr\bov(\theta)$, and that $\bl$ has weight $\heit(\bov(\theta)-\bov(\Theta))$.  In the case $|\bl|=0$ we have $\theta=\Theta$ and $w(\bl)=0$, so the result is trivial.  So we suppose $|\bl|>0$, and choose $i\in\bbz$ such that $\bl$ has a removable $i$-node.  We consider two cases.
\begin{description}
\item[$\theta_i(\bl)>0$\textnormal{:}]\indent\newline
Let $\bn = s_i(\bl)$.  Then $\bn$ has hub $s_i(\theta)$, and we have $w(\bn)=w(\bl)$, $\bov(s_i(\theta))=\bov(\theta)$ and $|\bn|=|\bl|+\theta_i(\bl)<|\bl|$, and the result follows by induction.
\item[$\theta_i(\bl)\ls0$\textnormal{:}]\indent\newline
Let $\bn$ be a multipartition obtained by removing a removable $i$-node from $\bl$, and let $\kappa$ be the hub of $\bn$.  Write $x=d_i(\theta)$ and $y = d_{i-1}(\theta)$.  Then $\theta_i=x-y$, and so by \cite[Lemma 3.6]{weight} we have $w(\bl)-w(\bn) = y-x+1$.  We also have $d_i(\kappa) = x-1$ and $d_{i-1}(\kappa)=y+1$, so that
\[v_j(\theta) = v_j(\kappa)+\delta_{jx}+\delta_{jy}-\delta_{j(x-1)}-\delta_{j(y+1)}.\]
Since $x\ls y$, this implies that $\bov(\theta) =\bov(\kappa)+\alpha_x+\alpha_{x+1}+\dots+\alpha_y$.  So $\bov(\kappa)\arr\bov(\theta)$, and $\heit(\bov(\theta)-\bov(\kappa))=w(\bl)-w(\bn)$; now the result follows by induction.
\end{description}
\end{pf}

\begin{eg}
As in the last example, suppose $r=4$ and $\boa = (1,2,1,0)$.  Let $\bl = ((1),(1^2),(1^2),(2))$.  The Young diagram of $\bl$ is
\[
\begin{picture}(30,45)
\put(0,20){\line(0,1){25}}
\put(15,30){\line(0,1){15}}
\put(0,30){\line(1,0){15}}
\put(0,45){\line(1,0){25}}
\put(5,19){$0$}
\put(5,34){$1$}
\put(20,34){$2$}
\end{picture}\qquad
\begin{picture}(30,45)
\put(0,5){\line(0,1){40}}
\put(15,15){\line(0,1){30}}
\put(0,15){\line(1,0){15}}
\put(0,30){\line(1,0){15}}
\put(0,45){\line(1,0){25}}
\put(5,4){$0$}
\put(5,19){$1$}
\put(5,34){$2$}
\put(20,34){$3$}
\end{picture}\qquad
\begin{picture}(30,45)
\put(0,5){\line(0,1){40}}
\put(15,15){\line(0,1){30}}
\put(0,15){\line(1,0){15}}
\put(0,30){\line(1,0){15}}
\put(0,45){\line(1,0){25}}
\put(2,4){$-1$}
\put(5,19){$0$}
\put(5,34){$1$}
\put(20,34){$2$}
\end{picture}\qquad
\begin{picture}(30,45)
\put(0,20){\line(0,1){25}}
\put(15,30){\line(0,1){15}}
\put(30,30){\line(0,1){15}}
\put(0,30){\line(1,0){30}}
\put(0,45){\line(1,0){40}}
\put(2,19){$-1$}
\put(5,34){$0$}
\put(20,34){$1$}
\put(35,34){$2$}
\end{picture}
\]
and this yields $\theta(\bl)=\kappa$, where
\[\kappa_i = \begin{cases}
-2 & (i=-1)\\
-1 & (i=0)\\
3 & (i=1)\\
-3 & (i=2)\\
-1 & (i=3)\\
0 & (\text{otherwise}).
\end{cases}\]
Note that $\kappa = s_1s_0\theta$, where $\theta$ is the hub from Example \ref{ex1}.  We have
\[\bov(\kappa)=\bov(\theta) = (-3,2,1,0,0)\]
and
\[\bov(\Theta) = (-2,1,0,1,0),\]
so $\bov(\kappa) = \bov(\Theta)+\alpha_1+\alpha_2$; and indeed, we may verify that the weight of $\bl$ is $2$.
\end{eg}

\sect{The case $e<\infty$}

For this section we fix $e\in\{2,3,\dots\}$ and $\boa\in(\zez)^r$.  As in Section \ref{sec3}, we answer our main question by considering negative hubs first, and then generalise using the $W_e$-action.  The results of Section \ref{sec3} play an important r\^ ole.

\subs{Each $W_e$-orbit contains one negative hub}\label{s41}

We begin by considering the orbits of the $W_e$-action on the set of hubs; we show how to identify the orbit containing a given hub, and prove that each orbit contains exactly one negative hub.

\begin{propn}\label{oneneghub}
Each $W_e$-orbit of hubs contains at least one negative hub.
\end{propn}

\begin{pf}
For any hub $\theta$, define
\[F(\theta) = \sum_{i,j\in\zez,\ i\neq j}(\theta_i+\theta_{i+1}+\dots+\theta_{j-1})^2.\]
We shall prove that if $\theta_k>0$ for some $k$, then $F(s_k(\theta))<F(\theta)$.  The proposition will then follow: given an orbit, a hub $\theta$ in that orbit for which the quantity $F(\theta)$ is minimised will necessarily be a negative hub.

So we assume that $\theta_k>0$, and let $\kappa = s_k(\theta)$.  Recall that this means that
\[\left(s_k(\theta)\right)_j = \theta_j + \theta_k\left(\delta_{(k-1)j}-2\delta_{kj}+\delta_{(k+1)j}\right)\]
for each $j$.  Putting $f_{ij}(\theta) = (\theta_i+\theta_{i+1}+\dots+\theta_{j-1})^2$, we have
\[f_{ij}(\theta) = f_{ij}(\kappa)\]
if the range $(i,i+1,\dots,j-1)$ contains all or none of $k\miu1,k,k\plu1$.  We also have
\begin{alignat*}3
f_{ik}(\theta) &= f_{i(k+1)}(\kappa),&\qquad f_{i(k+1)}(\theta)&=f_{ik}(\kappa)&\qquad&\text{for $i\neq k,k\plu1$}\\
f_{kj}(\theta) &= f_{(k+1)j}(\kappa),&\qquad f_{(k+1)j}(\theta)&=f_{kj}(\kappa)&\qquad&\text{for $j\neq k,k\plu1$}
\end{alignat*}
and
\[f_{k(k+1)}(\theta) = f_{k(k+1)}(\kappa).\]
So
\begin{align*}
F(\theta)-F(\kappa) &= f_{(k+1)k}(\theta)-f_{(k+1)k}(\kappa)\\
&=(\theta_{k+1}+\theta_{k+2}+\dots+\theta_{k-2}+\theta_{k-1})^2-\left((\theta_k+\theta_{k+1})+\theta_{k+2}+\dots+\theta_{k-2}+(\theta_{k-1}+\theta_k)\right)^2\\
&=(-r-\theta_k)^2 - (-r+\theta_k)^2\\
&= 4r\theta_k\\
&> 0,
\end{align*}
as required.
\end{pf}

To show that different negative hubs lie in different orbits, we introduce a set of invariants.  Given a hub $\theta$ and $i\in\zez$, define $t_i\in\bbz/er\bbz$ by
\[t_i(\theta) = ri+\theta_{i-1}+2\theta_{i-2}+\dots+(e-1)\theta_{i+1}.\]

The following simple observation will be very useful.

\begin{lemma}\label{useful}
For each $i\in\zez$ we have
\[t_{i}(\theta)-t_{i-1}(\theta) = -e\theta_i+er\bbz.\]
\end{lemma}

Now we show that the multiset $\{t_i(\theta)\mid i\in\zez\}$ is an invariant of the action of $W_e$ on hubs.  Then we show that if $\theta$ is negative then it is uniquely determined by this multiset, which implies that any orbit of hubs contains at most one negative hub.

\begin{propn}\label{tspres}
Suppose $\theta$ and $\kappa$ are hubs lying in the same $W_e$-orbit.  Then there is a permutation $\pi$ of $\zez$ such that
\[t_i(\kappa) = t_{\pi(i)}(\theta)\]
for all $i$.
\end{propn}

\begin{pf}
It suffices to consider the case where $\kappa = s_k(\theta)$ for some $k$.  In this case, it is easy to see that we have $t_i(\kappa)=t_i(\theta)$ unless $i=k-1$ or $k$, while
\begin{align*}
t_{k-1}(\kappa) &= rk-r + \theta_{k-2}+2\theta_{k-3}+\dots+(e-3)\theta_{k+2}+(e-2)(\theta_k+\theta_{k+1})+(e-1)(-\theta_k)\\
&= rk+\left(\textstyle\sum_{j\in\zez}\theta_j\right) + \theta_{k-2}+2\theta_{k-3}+\dots+(e-2)(\theta_{k+1})-\theta_k\\
&= rk+\theta_{k-1}+2\theta_{k-2}+\dots+(e-1)\theta_{k+1}\\
&= t_k(\theta);
\end{align*}
similarly we have $t_k(\kappa)=t_{k-1}(\theta)$, and the result follows.
\end{pf}

\begin{propn}\label{orbit}
Suppose $\theta$ and $\kappa$ are negative hubs, and that there is a permutation $\pi$ of $\zez$ such that
\[t_i(\kappa) = t_{\pi(i)}(\theta)\]
for all $i$.  Then $\theta=\kappa$.
\end{propn}

To prove Proposition \ref{orbit}, we consider cyclic ordering.  Given a positive integer $n$ and an $e$-tuple $t=(t_i)_{i\in\zez}$ of elements of $(\znz)^e$, we say that $t$ is \emph{in cyclic order modulo $n$} if there exist non-negative integers $(j_i)_{i\in\zez}$ summing to $n$ such that $t_i-t_{i-1} =  j_i+n\bbz$ for all $i$.

\begin{eg}
Suppose $e=4$ and $n=8$, and write $(t_i)_{i\in\bbz/4\bbz}$ as $(t_{\ol 0},t_{\ol 1},t_{\ol 2},t_{\ol 3})$.  Then $(0+8\bbz,0+8\bbz,3+8\bbz,3+8\bbz)$ and $(0+8\bbz,1+8\bbz,4+8\bbz,7+8\bbz)$ are in cyclic order modulo $8$, but $(0+8\bbz,3+8\bbz,0+8\bbz,3+8\bbz)$ and $(0+8\bbz,1+8\bbz,6+8\bbz,5+8\bbz)$ are not.
\end{eg}

\begin{pfof}{Proposition \ref{orbit}}
First we deal with the case where all the $t_i(\theta)$ (and hence all the $t_i(\kappa)$) are equal.  In this case, Lemma \ref{useful} implies that each $\theta_i$ is divisible by $r$; this means that $\theta_j$ equals $-r$ for some $j$, and $\theta_i=0$ for all $i\neq j$.  Similarly, we have $\kappa_k=-r$ for some $k$, with $\kappa_i=0$ for $i\neq k$.  But now the equality $t_j(\theta)= t_k(\kappa)$ gives $rj= rk$ (as elements of $\zerz$), so that $j=k$, and hence $\theta=\kappa$.

Now we turn to the case where the $t_i(\theta)$ are not all equal.  This means that each $\theta_i$ and each $\kappa_i$ lies in the range $\{0,-1,\dots,1-r\}$.  Since the integers $-e\theta_i$ are non-negative integers summing to $er$, Lemma \ref{useful} implies that the $e$-tuple $(t_i(\theta))_{i\in\zez}$ is in cyclic order modulo $er$.  Similarly, $(t_i(\kappa))_{i\in\zez}$ is in cyclic order modulo $er$.  This means that we may take the permutation $\pi$ to be a cyclic shift, i.e.\ there exists $k\in\zez$ such that
\[t_i(\kappa) = t_{i+k}(\theta)\]
for all $i$.  But now Lemma \ref{useful} gives
\[-e\kappa_i+er\bbz = t_{i}(\kappa)-t_{i-1}(\kappa) = t_{i+k}(\theta)-t_{i+k-1}(\theta) = -e\theta_{i+k}+er\bbz\]
for each $i$; since $0\gs\theta_i,\kappa_i\gs1-r$, we get
\[\kappa_i = \theta_{i+k}\]
for each $i$.  This in turn implies that
\begin{align*}
t_i(\kappa) &= ri + \theta_{i+k-1}+2\theta_{i+k-2}+\dots+(e-1)\theta_{i+k+1}\\
&=r(i+k) + \theta_{i+k-1}+2\theta_{i+k-2}+\dots+(e-1)\theta_{i+k+1}-rk\\
&= t_{i+k}(\theta)-rk\\
&=t_i(\kappa)-rk.
\end{align*}
We deduce that $k=\ol 0$, so $\kappa=\theta$.
\end{pfof}

\begin{eg}
Suppose $e=2$.  If $\theta$ is a hub, then we have
\begin{align*}
t_0(\theta) &= \theta_{\ol 1}+2r\bbz,\\
t_1(\theta) &= r+\theta_{\ol 0}+2r\bbz = -t_0(\theta).
\end{align*}
There are exactly $r+1$ different unordered pairs $(x,-x)$ for $x\in\bbz/2r\bbz$, and hence $r+1$ orbits of $W_2$ of the set of hubs.  And it is easy to see that there are exactly $r+1$ different negative hubs.
\end{eg}

\subs{The existence of a block with a given hub}

Now we turn to the problem of finding a necessary and sufficient condition for the existence of a block with a given hub.  We shall prove the following theorem, in which $\Theta$ is the hub $\Theta(e;\boa)$.

\begin{thm}\label{exi}
Suppose $\theta$ is a hub.  The following are equivalent.
\begin{enumerate}
\item
There exists a block with hub $\theta$.
\item
$\Theta\stackrel e\arr\theta$.
\item
$\left(\displaystyle\sum_{i\in\zez}i\theta_i\right)+e\bbz = \left(\displaystyle\sum_{i\in\zez}i\Theta_i\right)+e\bbz$.
\end{enumerate}
\end{thm}

First we observe that the second condition is invariant under the action of $W_e$.

\begin{lemma}\label{orbi}
Suppose $\theta$ and $\kappa$ are hubs lying in the same $W_e$-orbit.  Then $\theta\stackrel e\arr\kappa$.
\end{lemma}

\begin{pf}
It suffices to assume that $\kappa=s_l\theta$ for some $l\in\zez$.  Then (from Lemma \ref{hubaction}) we have $\kappa = \theta - \theta_i\alpha_i$.
\end{pf}

\begin{pfof}{Theorem \ref{exi}}
\begin{description}
\item[\textnormal{(1)$\Rightarrow$(2)}]
This is proved in exactly the same way as the implication (1)$\Rightarrow$(4) in Proposition \ref{tfae}.
\item[\textnormal{(2)$\Rightarrow$(3)}]
This is straightforward, since for any $j$ we have
\[\left(\sum_{i\in\zez}i(\alpha_j)_i\right)+e\bbz = \ol 0.\]
\item[\textnormal{(3)$\Rightarrow$(2)}]
We write $\zeta = \theta-\Theta$, and we must show that $\zeta$ is a linear combination of the $\alpha_i$, with integer coefficients.  Note that since $\Theta$ and $\theta$ are both hubs, we have
\[\sum_{i\in\zez}\zeta_i = 0,\]
and (assuming (3)) we have
\[\left(\sum_{i\in\zez}i\zeta_i\right) +e\bbz= \ol 0.\]
We define $m_{\ol 0}=0$ and
\[m_{\ol x} = \frac1e\left(\sum_{y=1}^xy(e-x)\zeta_{\ol y} + \sum_{y=x+1}^{e-1}x(e-y)\zeta_{\ol y}\right)\]
for $1\ls x\ls e-1$.  First we observe that each $m_{\ol x}$ is an integer, since we have
\[\sum_{y=1}^xy(e-x)\zeta_{\ol y} + \sum_{y=x+1}^{e-1}x(e-y)\zeta_{\ol y}\equiv\sum_{y=1}^{e-1}(-xy)\zeta_{\ol y}\equiv -x\sum_{y=0}^{e-1}y\zeta_{\ol y}\equiv0\pmod{e}.\]
Now we set $\zeta' = \sum_im_i\alpha_i$, and show that $\zeta'=\zeta$, i.e. $\zeta'_{\ol x}=\zeta_{\ol x}$ for all $x\in\{0,\dots,e-1\}$.  For $2\ls x\ls e-2$, we have
\begin{align*}
\zeta'_{\ol x} &= -m_{\ol x-1}+2m_{\ol x}-m_{\ol x+1}\\
&= \frac1e\Bigg(-\sum_{y=1}^{x-1}y(e-x+1)\zeta_{\ol y}+2\sum_{y=1}^{x-1}y(e-x)\zeta_{\ol y}-\sum_{y=1}^{x-1}y(e-x-1)\zeta_{\ol y}\\
&\phantom{=\frac1e\Bigg(}-\sum_{y=x+1}^{e-1}(x-1)(e-y)\zeta_{\ol y}+2\sum_{y=x+1}^{e-1}x(e-y)\zeta_{\ol y}-\sum_{y=x+1}^{e-1}(x+1)(e-y)\zeta_{\ol y}\\
&\phantom{=\frac1e\Bigg(}-(x-1)(e-x)\zeta_{\ol x}+2x(e-x)\zeta_{\ol x}-x(e-x-1)\zeta_{\ol x}\Bigg)\\
&= \zeta_{\ol x}.
\end{align*}
The cases $x=0,1,e-1$ are similar but simpler.
\item[\textnormal{(2)$\Rightarrow$(1)}]
Using Lemma \ref{orbi} and the transitivity and symmetry of $\stackrel e\arr$, we may replace $\theta$ with any hub in the same $W_e$-orbit.  In particular, since each orbit of hubs contains a negative hub, we may assume that $\theta$ is negative.

Assuming (2), let $(m_i)_{i\in\zez}$ be the unique $e$-tuple of non-negative integers such that $\theta=\Theta+\sum_im_i\alpha_i$ and some $m_k$ is zero.  Choose such a $k$, and let $\hat k$ be an integer such that $k=\ol{\hat k}$.  Now for each $i\in\zez$ define $\hat{\imath}$ to be the unique integer in the range $\{\hat k,\dots,\hat k+e-1\}$ such that $i=\ol{\hat{\imath}}$.  Define $\hat{\boa} = (\hat{a_1},\dots,\hat{a_r})\in\bbz^r$.  Now let $\hat{\Theta}$ be the $\infty$-hub $\Theta(\infty,\hat{\boa})$; that is, $\hat{\Theta}_{\hat{\imath}}=\Theta_i$ for each $i\in\zez$, and $\hat{\Theta}_j=0$ for all other $j$.  Define the $\infty$-hub $\hat{\theta}$ to be
\[\hat{\theta} = \hat{\Theta} + \sum_{i\in\zez}m_i\alpha_{\hat{\imath}}\]
Obviously $\hat{\Theta}\stackrel{\infty}\arr\hat{\theta}$, so by Proposition \ref{tfae} there is an $(\infty;\hat{\boa})$-block with hub $\hat{\theta}$.  If we take a multipartition $\bl$ lying in this block, then the $(e;\boa)$-block $B$ containing $\bl$ clearly has hub $\theta$.
\end{description}
\end{pfof}

\begin{eg}\label{ex2}
Suppose $r=4$, $e=3$ and $\boa=(\ol1,\ol2,\ol1,\ol0)$.  Then the hub $\Theta=\Theta(e;\boa)$ satisfies
\begin{alignat*}3
\Theta_{\ol0}&=-1,&\qquad\Theta_{\ol1}&=-2,&\qquad\Theta_{\ol2}&=-1.\\
\intertext{Let $\theta$ be the hub with}
\theta_{\ol0}&=-3,&\qquad\theta_{\ol1}&=-1,&\qquad\theta_{\ol2}&=0.
\end{alignat*}
We shall verify the conditions of Theorem \ref{exi} for $\theta$.  We have $\theta = \Theta+\alpha_{\ol1}+\alpha_{\ol2}$, so (2) holds.  For (3), we have
\[\left(\displaystyle\sum_{i\in\bbz/3\bbz}i\theta_i\right)+3\bbz = \ol2= \left(\displaystyle\sum_{i\in\bbz/3\bbz}i\Theta_i\right)+3\bbz.\]
For (1), we follow the proof of Theorem \ref{exi}.  Taking $k=\ol0$, we choose $\hat k=0$, so that $\hat\Theta$ and $\hat\theta$ are the hubs called $\Theta$ and $\theta$ in Example \ref{ex1}.  Taking the multipartition $\bm=(\odot,(1),(1),\odot)$ from that example, we see that the $(e;\boa)$-hub of $\bm$ is $\theta$.
\end{eg}

\subs{The weight of a core block}

In this section, we suppose that $\theta$ is a hub satisfying the conditions of Theorem \ref{exi}.  This means that there is a block, and hence a unique core block, with hub $\theta$.  We compute the weight of this core block; via Lemma \ref{listwt}, this enables us to identify all possible pairs $(\theta,w)$ such that there is a block of weight $w$.

Suppose throughout this section that $\theta$ is a hub satisfying the conditions of Theorem \ref{exi}.  We begin, as usual, with the case where $\theta$ is negative.  In this case, we let $k$, $\hat k$, $\hat\boa$ and $B$ be as in the proof of Theorem \ref{exi}.

\begin{lemma}\label{suffcore}
$B$ is a core block.
\end{lemma}

\begin{pf}
By construction we have $\hat\Theta_i=0$ unless $\hat k\ls i\ls \hat k+e-1$, and hence $\hat{\theta}_i=0$ unless $\hat k\ls i\ls \hat k+e$.  (For the fact that $\hat\theta_{\hat k-1}=0$, recall that $m_k=0$.)  So if we let $\bl$ be a multipartition in $B$ and calculate the beta-sets $B^1,\dots,B^r$ for $\bl$ using the multi-charge $\hat{\boa}$, then by Lemma \ref{recoverhub}(2) we must have
\[\{i\in\bbz\mid i\ls \hat k-1\}\subseteq B^j\subseteq\{i\in\bbz\mid i\ls \hat k+e-1\}\]
for each $j$.
Hence for $i\in\zez$ and $j\in\{1,\dots,r\}$ we have $\fb_{ij}^{\hat{\boa}}(\bl)\in\{\hat{\imath}-e,\hat{\imath}\}$, so $B$ is a core block.
\end{pf}

Now we give an expression for the weight of our core block $B$, analogous to the expression $\heit(\bov(\theta)-\bov(\Theta))$ of Theorem \ref{maininf}.  First we note the following lemma; this is a special case of a fact which is mentioned in \cite[\S3.4]{core}.

\begin{lemma}\label{wtsame}
The $(e;\boa)$-weight of $B$ equals the $(\infty;\hat\boa)$-weight of $B$.
\end{lemma}

\begin{pf}
Let $\bl$ be a multipartition in $B$.  We claim that the $(\infty;\hat\boa)$-residue of any node of $\bl$ lies in the set $\{\hat k+1,\dots, \hat k+e-1\}$.  To prove the claim, we assume $\bl\neq\varnothing$ and let $m$ be the maximum $(\infty;\hat\boa)$-residue of any node.  If the node $(i,j)_l$ is a node of $\bl$ with residue $m$, then the node $(i,j+1)_l$ is an addable node of $\bl$ with residue $m+1$.  Since there are no removable nodes of residue $m+1$, we must have $\hat\theta_{m+1}<0$.  But from the proof of Lemma \ref{suffcore} we have $\hat\theta_i=0$ unless $i\ls \hat k+e$, so we must have $m\ls\hat k+e-1$.  Similarly the minimal $(\infty;\hat\boa)$-residue of any node of $\bl$ is at least $\hat k+1$, and the claim is proved.

We see that $\bl$ has no nodes of $(e;\boa)$-residue $k$, and that for $i\neq k$ a node of $\bl$ has $(e;\boa)$-residue $i$ if and only if it has $(\infty;\hat\boa)$-residue $\hat\imath$.  Now it is clear from the formula for weight that the $(e;\boa)$-weight and the $(\infty;\hat\boa)$-weight coincide.
\end{pf}

Recall from \S\ref{s41} the definition of $t_i(\theta),t_i(\Theta)\in\zerz$.  By Lemma \ref{useful} we see that all the $t_i(\theta)$ are congruent modulo $e$, meaning that for any $i,j$ there is an integer $x$ such that $t_i(\theta)-t_j(\theta) = ex+er\bbz$.  Similarly, all the $t_i(\Theta)$ are congruent modulo $e$.  In fact, our assumption that $\Theta\stackrel e\arr\theta$ implies that the $t_i(\Theta)$ are congruent to the $t_i(\theta)$ modulo $e$.  To see this, recall that we have
\[\theta = \Theta+\sum_{i\in\zez}m_i\alpha_i\]
with $m_k=0$; comparing this with the definitions of $t_k(\theta)$ and $t_k(\Theta)$ yields
\begin{equation}
t_k(\theta)-t_k(\Theta) = em_{k+1}+er\bbz,\tag*{($\ast$)}
\end{equation}
so $t_k(\theta)$ and $t_k(\Theta)$ are congruent modulo $e$.

Let $h$ be the unique integer in the set $\{0,\dots,e-1\}$ such that $\ol{t_i(\theta)}=\ol{t_i(\Theta)}=\ol h$ for each $i$.  Then both the expressions $h-t_i(\Theta)$ and $h-t_i(\theta)$ are divisible by $e$; given $j\in\zrz$, we write
\begin{align*}
v_j(\Theta) &= \left|\left\{i\in\zez\ \left|\ h-t_i(\Theta) = ej\right.\right\}\right|,\\
v_j(\theta) &= \left|\left\{i\in\zez\ \left|\ h-t_i(\theta) = ej\right.\right\}\right|.
\end{align*}
Now define $\bov(\theta)$ to be the element of $\bbz^{\zrz}$ with coordinates $v_j(\theta)$.  Define $\bov(\Theta)$ analogously.

\begin{propn}\label{compwt}
Suppose $\theta$ is a negative hub such that $\Theta\stackrel e\arr\theta$, and define $\bov(\theta),\bov(\Theta)$ as above.  Then $\bov(\Theta)\stackrel r\arr\bov(\theta)$, and the block $B$ has weight $\heit(\bov(\theta)-\bov(\Theta))$.
\end{propn}

\begin{pf}
Let $\bl$ be a multipartition in $B$.  By Lemma \ref{wtsame}, the $(e;\boa)$-weight of $\bl$ equals its $(\infty;\hat{\boa})$-weight, and by Theorem \ref{mainneg}, this is $\heit(\bov(\hat{\theta})-\bov(\hat{\Theta}))$.  So we must show that $\bov(\Theta)\stackrel r\arr\bov(\theta)$ and $\heit(\bov(\theta)-\bov(\Theta))=\heit(\bov(\hat\theta)-\bov(\hat\Theta))$.

Recall the definition of the integers $d_l(\hat\theta)$ for $l\in\bbz$.  By construction, we have
\[d_l(\hat\theta) = \begin{cases}
r & (l\ls \hat k-1)\\
r+\hat\theta_{\hat k}+\hat{\theta}_{\hat k+1}+\dots+\hat\theta_l & (l=\hat\imath\text{, some }i\in\zez)\\
0 & (l\gs \hat k+e),
\end{cases}\]
with a similar expression for $\hat\Theta$.  Recalling the construction of $\hat\Theta$ and $\hat\theta$, we see that For $i\in\zez$, we have
\begin{equation}
\hat\Theta_{\hat{\imath}} = \Theta_i,\qquad \hat\theta_{\hat{\imath}} = \begin{cases}
\Theta_k-m_{k+1} & (i=k)\\
\theta_i & (i\neq k).
\end{cases}\tag*{($\dagger$)}
\end{equation}

Let $h$ be as above, and let $\hslash$ be the element of $\zrz$ such that $h-t_k(\Theta)=e\hslash$.  Now take $1\ls x\ls r$, and write $\ul x = x-\Theta_k+\hslash$.  Since $0\ls d_l(\Theta)\ls r$ for each $l$, we have $d_l(\Theta)=x$ if and only if $d_l(\Theta)\equiv x\ppmod r$.  Hence
\begin{align*}
v_{\ul{x}}(\Theta) &= \left|\left\{i\in\zez\ \left|\ h-t_i(\Theta) = e\ul x\right.\right\}\right|\\
&= \left|\left\{i\in\zez\ \left|\ t_k(\Theta)-t_i(\Theta) = e\ul x+t_k(\Theta)-h\right.\right\}\right|\\
&= \left|\left\{i\in\zez\ \left|\ e(\Theta_{k+1}+\dots+\Theta_i)+er\bbz = e\ul x+t_k(\Theta)-h\right.\right\}\right|\tag*{by Lemma \ref{useful}}\\
&= \left|\left\{i\in\zez\ \left|\ \Theta_{k+1}+\dots+\Theta_i+r\bbz = \ul x-\hslash\right.\right\}\right|\\
&= \left|\left\{i\in\zez\ \left|\ \hat\Theta_{\hat k}+\dots+\hat\Theta_{\hat\imath}+r\bbz = \ul x-\hslash+\Theta_k\right.\right\}\right|\tag*{by ($\dagger$)}\\
&= \left|\left\{i\in\zez\ \left|\ d_{\hat\imath}(\hat\Theta)\equiv x\pmod r\right.\right\}\right|\\
&= \left|\left\{l\in\bbz\ \left|\ d_l(\hat\Theta)=x\right.\right\}\right|\\
&=v_x(\hat\Theta).
\end{align*}
For $\theta$, we have
\begin{align*}
v_{\ul{x}}(\theta) &= \left|\left\{i\in\zez\ \left|\ h-t_i(\theta) = e\ul x\right.\right\}\right|\\
&= \left|\left\{i\in\zez\ \left|\ t_k(\theta)-t_i(\theta)-em_{k+1}=e\ul x+t_k(\theta)-em_{k+1}-h\right.\right\}\right|\\
&= \left|\left\{i\in\zez\ \left|\ e(\theta_{k+1}+\dots+\theta_i-m_{k+1}) +er\bbz = e\ul x+t_k(\Theta)-h\right.\right\}\right|\tag*{by Lemma \ref{useful} and ($\ast$)}\\
&= \left|\left\{i\in\zez\ \left|\ \theta_{k+1}+\dots+\theta_i-m_{k+1}+r\bbz = \ul x-\hslash\right.\right\}\right|\\
&= \left|\left\{i\in\zez\ \left|\ \hat\theta_{\hat k}+\dots+\hat\theta_{\hat\imath}+r\bbz = \ul x-\hslash+\Theta_k\right.\right\}\right|\tag*{by ($\dagger$)}\\
&= \left|\left\{i\in\zez\ \left|\ d_{\hat\imath}(\hat\theta)\equiv x\pmod r\right.\right\}\right|\\
&= \left|\left\{l\in\bbz\ \left|\ d_l(\hat\theta) = x\right.\right\}\right|\\
&=v_x(\hat\theta).
\end{align*}
Also, since
\[\sum_{x=0}^rv_x(\hat\Theta) = 0\qquad\text{and}\qquad\sum_{j\in\zrz}v_j(\Theta)=e,\]
we have
\begin{align*}
v_{\hslash-\Theta_k}(\Theta) &= v_0(\hat\Theta)+v_r(\hat\Theta)+e;\\
\intertext{similarly}
v_{\hslash-\Theta_k}(\theta) &= v_0(\hat\theta)+v_r(\hat\theta)+e.
\end{align*}
Now we can complete the proof: the statement $\bov(\hat\Theta)\stackrel{0,r}\arr\bov(\hat\theta)$ says that there are non-negative integers $n_1,\dots,n_{r-1}$ such that
\begin{align*}
\bov(\hat\theta) &= \bov(\hat\Theta)+\sum_{x=1}^{r-1}n_x\alpha_x,\\
\intertext{and this combined with the above expressions gives}
\bov(\theta) &= \bov(\Theta)+\sum_{x=1}^{r-1}n_i\alpha_{\ul x}.
\end{align*}
This yields an expression for $v(\Theta)-v(\theta)$ as a linear combination of the $(\alpha_i)_{i\in\zrz}$ with all coefficients non-negative integers and at least one coefficient (namely, the coefficient of $\alpha_{\hslash-\Theta_k}$) equal to zero.  Hence we have $\bov(\Theta)\stackrel r\arr\bov(\theta)$, and
\[\heit(\bov(\theta)-\bov(\Theta)) = \sum_{x=1}^{r-1}n_x = \heit(\bov(\hat\theta)-\bov(\hat\Theta)),\]
as required.
\end{pf}

\begin{eg}
Retain the notation from Example \ref{ex2}.  We may calculate
\begin{alignat*}3
t_{\ol0}(\Theta) &= 7+12\bbz,&\qquad t_{\ol1}(\Theta) &= 1+12\bbz,&\qquad t_{\ol2}(\Theta) &=4+12\bbz,\\
t_{\ol0}(\theta) &= 10+12\bbz,&\qquad t_{\ol1}(\theta) &= 1+12\bbz,&\qquad t_{\ol2}(\theta) &=1+12\bbz.
\end{alignat*}
So $h=1$, and
\begin{alignat*}4
v_{0+4\bbz}(\Theta) &= 1,&\qquad v_{1+4\bbz}(\Theta) &= 0,&\qquad v_{2+4\bbz}(\Theta) &= 1,&\qquad v_{3+4\bbz}(\Theta) &= 1,\\
v_{0+4\bbz}(\theta) &= 2,&\qquad v_{1+4\bbz}(\theta) &= 1,&\qquad v_{2+4\bbz}(\theta) &= 0,&\qquad v_{3+4\bbz}(\theta) &= 0.
\end{alignat*}
We see that $\bov(\theta)=\bov(\Theta)+\alpha_{0+4\bbz}+\alpha_{1+4\bbz}$, so $\bov(\Theta)\stackrel4\arr\bov(\theta)$ and $\heit(\bov(\theta)-\bov(\Theta))=2$.  And indeed, the $(e;\boa)$-weight of $\bm$ is $2$.
\end{eg}

Now we can generalise to the case where $\theta$ is not necessarily negative, and give our main theorem about the weight of a core block, in the case where $e$ is finite.

\begin{thm}\label{finmain}
Suppose $e<\infty$ and $\theta$ is a hub satisfying the conditions of Theorem \ref{exi}, and define $\bov(\theta)$ as above.  Then $\bov(\Theta)\stackrel{r}\arr\bov(\theta)$, and the weight of the core block with hub $\theta$ is $\heit(\bov(\theta)-\bov(\Theta))$.
\end{thm}

\begin{pf}
In the case where $\theta$ is negative, the theorem follows from Lemma \ref{suffcore} and Proposition \ref{compwt}.  To extend to arbitrary hubs, we note that the following are preserved under the $W_e$-action on hubs:
\begin{itemize}
\item
the existence or not of a block with hub $\theta$ (Corollary \ref{existorb}(\ref{eo1}));
\item
the weight of the core block with hub $\theta$, if one exists (Corollary \ref{existorb}(\ref{eo2}));
\item
the multiset $\{t_i(\theta)\mid i\in\zez\}$ (Proposition \ref{tspres}), and hence $\bov(\theta)$.
\end{itemize}
Hence by Proposition \ref{oneneghub} the result holds for all hubs.
\end{pf}

\end{document}